\documentclass[12pt,dvips]{amsart}
\usepackage{amsfonts, amssymb, latexsym, epsfig}

\newcommand{\End}{\operatorname{End}}

\newtheorem{theorem}{Theorem}
\newtheorem{corollary}{Corollary}

\begin{document}

\title{Interval Pattern Avoidance for Arbitrary Root Systems}
\subjclass[2000]{14M15; 05E15}
\author{Alexander Woo}
\address{Department of Mathematics, Mathematical Sciences Building, One Shields Ave., 
University of California, Davis, CA, 95616, USA}
\email{awoo@math.ucdavis.edu}

\date{\today}
\begin{abstract}

We extend the idea of interval pattern avoidance defined by Yong
and the author for $S_n$ to arbitrary Weyl groups using the definition
of pattern avoidance due to Billey and Braden, and Billey and
Postnikov.  We show that, as previously shown by Yong and the
author for $GL_n$, interval pattern avoidance is a universal tool for
characterizing which Schubert varieties have certain local properties,
and where these local properties hold.

\end{abstract}
\maketitle

\section{Introduction}

The purpose of this brief note is to extend the notions of interval
pattern embedding and avoidance introduced by Yong and the
author in type A \cite{WYII} to Schubert varieties of arbitrary Lie
type.  This extension is the natural common generalization of the
definition in type A and the definition of pattern avoidance coming
from root subsystems, as introduced combinatorially by Billey and
Postnikov \cite{BilPos} and explained geometrically via the pattern
map by Billey and Braden \cite{BilBra}.  (In type A, the pattern map
was also implicit in work of N. Bergeron and Sottile \cite{BerSot}.)

The main reason for our definition of interval pattern avoidance is
that it gives a universal tool for describing local properties on
Schubert varieties, in the sense that the set of points on all
Schubert varieties satisfying any given local property (except for
dimension) has a characterization using only interval pattern
avoidance.  The main example of such a property for which results are
known is smoothness.  The Schubert varieties which are smooth
everywhere is characterized by ordinary pattern avoidance
\cite{LakSand, Bil, BilPos}.  The locus of singular points in any
Schubert variety of type A was described independently in
\cite{BilWar, Cor, KLR, Man}; this description can be easily
reformulated in terms of interval pattern embeddings
\cite[Thm. 6.1]{WYII}, and interval pattern embeddings should provide
the appropriate language for a similar description in general.
However, the Schubert varieties which are everywhere Gorenstein
\cite{WYI} or everywhere factorial \cite{BousBut} cannot be
characterized by ordinary pattern avoidance.  Nevertheless, interval
pattern avoidance suffices, not only for these properties, but for any
local property preserved under products with affine space.

This universality is demonstrated by showing that interval pattern
embeddings give an isomorphism of slices of different Schubert
varieties.  This isomorphism is proven using the pattern map of Billey
and Braden; when written in coordinates, the proof becomes essentially
the same as the one previously given for type A.  However, this new
proof shows the isomorphism extends to the Richardson varieties which
are the closures of the slices, which is a new result even in type A.

I thank Frank Sottile and Mark Haiman, who independently suggested
that there should be some geometric explanation for the forerunner of
interval pattern avoidance found in the characterization of Gorenstein
Schubert varieties.  In addition, I thank Alexander Yong and Sara
Billey for helpful discussions, and Vic Reiner and Ezra Miller for
organizational advice and for the name ``interval pattern avoidance''.
Research was supported in part by NSF VIGRE grant No. DMS--0135345.

\section{Pattern Avoidance and Interval Pattern Avoidance}

Fix a ground field $\Bbbk$.  Let $G$ be a connected semisimple
linear algebraic group over $\Bbbk$, $B$ a fixed Borel subgroup,
and $T\subseteq B\subseteq G$ a maximal torus.  Let $\Phi$ be the
roots of $G$ under the action of $T$, $\Phi_+$ and $\Phi_-$ the
positive and negative roots corresponding to our choice of Borel
subgroup, and $\Delta_+$ the set of simple positive roots.  Let $V$ be
the inner product space spanned by the root lattice.  The Weyl group
$W$ of $G$ is the group generated by the reflections
$\{s_\alpha\in\End(V)\mid \alpha\in \Phi\}$, where
$s_\alpha(v)=v-2(v,\alpha)/(\alpha,\alpha)\alpha$.  The {\bf length}
$\ell(w)$ of an element in a Weyl group is the minimal length of any
expression $w=s_{\beta_1}s_{\beta_2}\cdots s_{\beta_k}$, where each
$\beta_j$ is a simple root.  The Weyl group can also be recovered from
$G$ as the group $N(T)/T$.  Pattern avoidance depends not only on the
abstract Weyl group but also on the root system it acts on; to
emphasize this, in the remainder of the paper, we denote a Weyl group
by the triple $(W,\Phi,V)$.

The variety $G/B$ is known as the {\bf flag variety}.  The group $G$
acts on $G/B$ via multiplication on the left.  To each element of
$u\in (W,\Phi,V)$ (here considered as $N(T)/T$) we can associate the
$T$-fixed point $e_u:=uB$, and these are all the $T$-fixed points of
$G/B$.  There is a {\bf Bruhat decomposition} of $G/B$ into {\bf
Schubert cells} $X^\circ_w:=Be_wB/B$, one for each $w\in W$, and the
{\bf Schubert variety} $X_w$ is the closure of the Schubert cell
$X^\circ_w$.  There is also a decomposition of $G/B$ into {\bf
opposite Schubert cells} $\Omega^\circ_w=B_-e_wB/B$, where $B_-$ is
the Borel subgroup opposite to $B$; the closure of the opposite
Schubert cell $\Omega^\circ_w$ is called an {\bf opposite Schubert
variety} and is denoted $\Omega_w$.  The {\bf Richardson variety}
$X^u_v$ is the intersection of $\Omega_u$ and $X_v$; Richardson showed
that it is reduced and irreducible (when nonempty) \cite{Rich}.  The
dimension of $X_w$ and the codimension in $G/B$ of $\Omega_w$ are both
$\ell(w)$.  The dimension of $X^u_v$ is $\ell(v)-\ell(u)$.

The Schubert variety $X_w$ is a union of Schubert cells.  We define
the {\bf Bruhat order} on $(W,\Phi,V)$ by declaring that $u\leq v$ if
$X^\circ_u\subseteq X_v$.  Alternatively, Bruhat order can be defined
combinatorially by declaring it to be the reflexive transitive closure
of the relation $\prec$ under which $u\prec v$ if both $u=s_\alpha v$
for some $\alpha\in\Phi$ and $\ell(u)<\ell(v)$.  This combinatorial
definition has a geometric explanation; when $u$ and $v$ are so
related, the curve $\overline{U_\alpha\cdot e_v}$ is a $\mathbb{P}^1$
inside $X_v$ connecting $e_u$ and $e_v$.  Here, $U_\alpha$ is the root
subgroup of $B$ corresponding to the root $\alpha$.  The Richardson
variety $X^u_v$ is nonempty whenever $u\leq v$.

Now we recall the definitions of pattern embeddings and pattern
avoidance found in \cite{BilBra, BilPos}.  Let
$(W^\prime,\Phi^\prime,V^\prime)$ and $(W,\Phi,V)$ be Weyl groups.
A {\bf subsystem embedding} $i$ of $(W^\prime,\Phi^\prime,V^\prime)$ into
$(W,\Phi,V)$ is an embedding of $V^\prime$ as a subspace of $V$ so that
$\Phi^\prime\cong\Phi\cap i(V^\prime)$; this induces an embedding of
$W^\prime$ into $W$ as the subgroup generated by the reflections
$\{s_\alpha\mid \alpha\in i(\Phi^\prime)\}$.

Define the {\bf flattening map} $\phi_i$ from $(W,\Phi,V)$ to
$(W^\prime,\Phi^\prime,V^\prime)$ as follows.  An element
$w\in(W,\Phi,V)$ is uniquely determined by its inversion set
$I(w)=\Phi_+\cap w(\Phi_-)$.  Therefore we can define $\phi_i(w)$ as
the element of $(W^\prime,\Phi^\prime,V^\prime)$ whose inversion set
is $i^{-1}(I(w)\cap i(\Phi^{\prime+}))$.  Then {\bf $i$ (pattern)
embeds $v\in(W^\prime,\Phi^\prime,V^\prime)$ in $w\in(W,\Phi,V)$} if
$\phi_i(w)=v$.  The Weyl group element $w$ is said to {\bf (pattern)
avoid} $v$ if $\phi_i(w)\neq v$ for every embedding $i$ of
$(W^\prime,\Phi^\prime,V^\prime)$ into $(W,\Phi,V)$.

Our definition of interval pattern avoidance is now as follows.  Let
$u\leq v\in (W^\prime,\Phi^\prime,V^\prime)$ and $x\leq w\in
(W,\Phi,V)$, where $\leq$ denotes the Bruhat order.  Let $i$ be a
subsystem embedding of $(W^\prime,\Phi^\prime,V^\prime)$ into
$(W,\Phi,V)$.  We say {\bf $i$ (interval pattern) embeds $[u,v]$ in
$[x,w]$} if the following three conditions are all satisfied.
\begin{enumerate}
\item $\phi_i(w)=v$ and $\phi_i(x)=u$.
\item $x$ and $w$ are in the same right $i(W^\prime)$ coset.
\item $[u,v]$ and $[x,w]$ are isomorphic as intervals in Bruhat order.
\end{enumerate}

The third condition implies in particular that
$\ell(v)-\ell(u)=\ell(w)-\ell(x)$.  This equality in lengths is
actually sufficient to imply the third condition, given the first two;
a combinatorial proof of this fact is possible, but the geometry below
also shows it.

Note that the first two conditions imply that $x=i(uv^{-1})w$.  Since
$x$ is determined by $u$, $v$, $w$, and $i$, we will say that {\bf $w$
(interval pattern) avoids $[u,v]$} if, for every subsystem embedding
$i$ of $(W^\prime,\Phi^\prime,V^\prime)$ into $(W,\Phi,V)$, $[u,v]$
does not embed in $[i(uv^{-1})w,w]$.

\section{Main Theorem and Corollary}

Our main theorem can now be stated as follows.

\begin{theorem}
Suppose there is some subsystem embedding $i$ which embeds $[u,v]$ in
$[x,w]$.  Then the Richardson varieties $R^u_v$ and $R^x_w$ are
isomorphic.  This isomorphism sends $\Omega^\circ_\sigma\cap
X^\circ_\tau$ to $\Omega^\circ_{\phi_i(\sigma)} \cap
X^\circ_{\phi_i(\tau)}$ for every $\sigma,\tau\in[x,w]$.
\end{theorem}

The main application of this theorem we have in mind is to the study
of singularities of Schubert varieties.  Call a local property
$\mathcal{P}$ {\bf semicontinuously stable} if it is preserved under
products with affine space, and the $\mathcal{P}$-locus on any
Schubert variety is closed.  Examples include being singular, being
non-Gorenstein, having multiplicity greater than some fixed number
$k$, or having a particular coefficient of the Kazhdan-Luzstig
polynomial being greater than a fixed number $k$.  Now define a poset
on the set of all intervals in all Weyl groups (where, as throughout,
the root system is considered part of the data of the Weyl group) by
taking the reflexive transitive closure of the following two
relations.

\begin{enumerate}
\item $[u,v]\prec[x,w]$ if there is some embedding of $[u,v]$ into
  $[x,w]$.
\item $[u,v]\prec[u^\prime,v]$ if $u\leq u^\prime$
\end{enumerate}

Now we can state our corollary.

\begin{corollary}
Let $\mathcal{P}$ be a semicontinuously stable property.  Then the set
of intervals such that $\{[u,v]\mid \mathcal{P} \mbox{ holds at
}e_u\mbox{ on }X_v\}$ is an upper order ideal on the aforementioned
poset.  The set $\{w\mid \mathcal{P} \mbox{ holds on no points of }
X_w\}$ is the set of $w$ avoiding some list of intervals $[u,v]$.
\end{corollary}

Notice that this corollary holds separately for different ground
fields, in that the order ideal for the same property may depend on
$\Bbbk$.  The list of intervals to be avoided may be infinite, but we
hope that for any particular property it has a nice form.

\begin{proof}
The point $e_u$ has a neighborhood $u\cdot\Omega^\circ_{\mathrm{id}}$
in $G/B$, so $u\cdot\Omega^\circ_{\mathrm{id}}\cap X_v$ is a
neighborhood of $e_u$ on $X_v$.  This neighborhood is isomorphic to
$(\Omega^\circ_u\cap X_v) \times \mathbb{A}^{\ell(u)}$ \cite[Lemma
A.4]{KazLus}.  Therefore, any semicontinuously stable property
$\mathcal{P}$ depends only on $\Omega^\circ_u\cap X_v$, which is
commonly called the {\bf slice} of $X_v$ at $e_u$.  Our theorem now
shows that $\mathcal{P}$ is preserved under going up in our poset by
the first type of generating relation, since $\Omega^\circ_u\cap X_v$
is isomorphic to $\Omega^\circ_x\cap X_w$.

As for the second type of generating relation, we can by induction on
Bruhat order reduce to the case where $u^\prime=s_\alpha u$.  In that
case, $U_\alpha\cdot e_u$ is a curve in $X_v$ all of whose
points have neighborhoods isomorphic to the neighborhood at $e_u$
(since $X_v$ has a $B$-action).  The closure of $U_\alpha\cdot e_u$
includes the additional point $e_u^\prime$.  Since the set at which
$\mathcal{P}$ holds is closed, $\mathcal{P}$ is also preserved going
up by the second type of generating relation.

The last statement follows by taking a generating set for the order ideal.
\end{proof}

We also have the following corollary about Kazhdan-Luzstig
polynomials, generalizing a lemma of Polo \cite[Lemma 2.6]{Polo}.
(See also \cite[Thm. 6]{BilBra}.)

\begin{corollary}
Suppose a subsystem embedding embeds $[u,v]$ into $[x,w]$.  Then the
Kazhdan-Luzstig polynomials $P_{u,v}(q)$ and $P_{x,w}(q)$ are equal.
\end{corollary}

It is conjectured that $P_{u,v}(q)=P_{x,w}(q)$ whenever $[u,v]$ and
$[x,w]$ are isomorphic as intervals, and this theorem confirms a very
special case of this conjecture.  Kazhdan-Luzstig polynomials and this
conjecture are discussed with further references in \cite{BjoBre}.

\section{The Pattern Map}

To prove the theorem, we use the geometric pattern map introduced by
Billey and Braden \cite{BilBra}.  Let $T_0$ be a one parameter
subgroup of $T$ which is generic among subgroups satisfying
$\alpha(T_0)=1$ for every $\alpha\in i(\Phi^\prime)$.  (Recall that
roots are actually irreducible representations of $T$, which are maps
from $T$ to $\Bbbk^\times$.)  Let $G^\prime$ be the centralizer
$Z_G(T_0)$ of $T_0$.  The Weyl group and roots of $G^\prime$ are then
$i(W^\prime)$ and $i(\Phi^\prime)$.  In $G^\prime$ we fix the Borel
subgroup $B^\prime=G^\prime\cap B$.

Now Billey and Braden define a map $\psi: (G/B)^{T_0} \rightarrow
(G^\prime/B^\prime)$ as follows.  There is a bijection between points
of $G/B$ and Borel subgroups of $G$ given by associating to the coset
$gB$ the Borel subgroup $gBg^{-1}$.  Now define $\psi(gB)$ to be the
point in $G^\prime/B^\prime$ associated with the Borel subgroup
$gBg^{-1}\cap G^\prime$.  This is a Borel subgroup of $G^\prime$
whenever $gB$ is fixed by $T_0$ \cite[Thm. 6.4.7]{Spr}.  Billey and
Braden prove the following theorem.

\begin{theorem}\cite[Thm. 10]{BilBra}
\begin{enumerate}
\item The map $\psi$ restricts to an isomorphism on each connected
component of $(G/B)^{T_0}$.
\item For any $w\in (W,\Phi,V)$, the restriction of $\psi$ is an
isomorphism between $X^\circ_w\cap(G/B)^{T_0}$ and
$X^\circ_{\phi_i(w)}$ taking $e_w$ to $e_{\phi_i(w)}$.
\end{enumerate}
\end{theorem}

Their proof of part 2 also shows that $\psi$ restricts to an
isomorphism between $\Omega^\circ_w \cap (G/B)^{T_0}$ and
$\Omega^\circ_{\phi_i(w)}$.

As remarked by Billey and Braden \cite{BilBra} (see also
\cite[Prop. 4.2]{LenRobSot}), this geometric pattern map explains why
ordinary pattern avoidance characterizes singular Schubert varieties.
Given any one parameter torus $T_0\cong\Bbbk^\times$ acting on a
Schubert variety, if the $T_0$ fixed locus is singular, the entire
Schubert variety must be singular.  If there is a pattern embedding
$i$ of $v$ into $w$, and $X_v$ is singular, then $\psi$ gives an
isomorphism between $X_w\cap (G/B)^{T_0}$ and $X_{v}$, showing that
$X_w$ is singular.

\section{Proof of Main Theorem}

We will show that $X^u_v$ and $X^x_w$ are isomorphic under the map
$\psi$.  First we show that $e_w$ and $e_x$ are in the same connected
component of $(G/B)^{T_0}$, as follows.  Since $x$ and $w$ are in the
same right $W^\prime$ coset, we can successively multiply $x$ on the
left by reflections $s_\alpha$, with $\alpha\in\Phi^\prime$, to get
$w$.  If $\sigma=s_\alpha\tau$ for some $\alpha\in\Phi^\prime$, then
$e_\sigma$ and $e_\tau$ are in the same connected component of
$(G/B)^{T_0}$ since the points are connected by the Schubert curve
$\overline{U_\alpha\cdot e_\sigma}$ (assuming $\sigma\geq\tau$), and
$U_\alpha$ is $T_0$ fixed as $\alpha\in\Phi^\prime$.

Now, by part 2 of the theorem, $X^\circ_w\cap(G/B)^{T_0}$ and
$\Omega^\circ_x\cap(G/B)^{T_0}$ are connected, so, given that $e_w$
and $e_x$ are in the same connected component of $(G/B)^{T_0}$,
$(X^\circ_w\cup\Omega^\circ_x)\cap(G/B)^{T_0}$ is contained in a
single connected component of $(G/B)^{T_0}$.  Therefore, $\psi$ is an
isomorphism when restricted to
$X^\circ_w\cap\Omega^\circ_x\cap(G/B)^{T_0}$.

We show that the image of $X^\circ_w\cap\Omega^\circ_x\cap(G/B)^{T_0}$
is $X^\circ_v\cap\Omega^\circ_u$.  If $p\in
X^\circ_v\cap\Omega^\circ_u$, then $p=\psi(p_1)$ for some $p_1\in
X^\circ_w\cap (G/B)^{T_0}$, and $p=\psi(p_2)$ for some
$p_2\in\Omega^\circ_x\cap (G/B)^{T_0}$.  Since $X^\circ_w\cap
(G/B)^{T_0}$ and $\Omega^\circ_x\cap (G/B)^{T_0}$ lie in the same
connected component of $(G/B)^{T_0}$, by part 1 of the theorem,
$p_1=p_2\in X^\circ_w\cap\Omega^\circ_x\cap(G/B)^{T_0}$.  Conversely,
for $p\in X^\circ_w\cap\Omega^\circ_x\cap(G/B)^{T_0}$, $\psi(p)\in
\psi(X^\circ_w\cap(G/B)^{T_0})\cap\psi(\Omega^\circ_x\cap(G/B)^{T_0})=X^\circ_v\cap\Omega^\circ_u$.

In particular, $X^\circ_w\cap\Omega^\circ_x\cap(G/B)^{T_0}$ has
dimension $\ell(v)-\ell(u)$.  Since we have a pattern embedding from
$[u,v]$ to $[x,w]$, $\ell(v)-\ell(u)=\ell(w)-\ell(x)$, so the
dimension of $X^\circ_w\cap\Omega^\circ_x\cap(G/B)^{T_0}$ is the same
as the dimension of $X^\circ_w\cap\Omega^\circ_x$.  As the latter is
known to be irreducible \cite{Rich} (or \cite[Prop.  1.3.2]{Bri}),
$X^\circ_w\cap\Omega^\circ_x$, and therefore its closure $X^x_w$, must
be pointwise $T_0$-fixed.

Since $X^x_w$ is connected and pointwise $T_0$-fixed, it must be
isomorphic to its image under $\psi$.  This image is the closure of
$X^\circ_v\cap\Omega^\circ_u$, which is $X^u_v$.

The calculation of the image can be repeated for every
$\sigma,\tau\in[x,w]$, proving the second statement.

\end{document}